\documentclass[11pt]{amsart}
\usepackage{amssymb,amscd}

\newtheorem{theorem}{Theorem}[section]

\newtheorem{proposition}[theorem]{Proposition}

\newtheorem{lemma}[theorem]{Lemma}

\theoremstyle{remark}
\newtheorem{remark}[theorem]{Remark}
\newtheorem{definition}[theorem]{Definition}

\newcommand\be{\begin{equation}\label}
\newcommand\ee{\end{equation}}
\newcommand\M{\mathcal{M}}

\newcommand{\N}{\mathcal{N}}
\newcommand{\R}{\mathbb{R}}

\newcommand{\QQ}{\mathbb{Q}}

\newcommand\lie[1]{\mathfrak{#1}}
\renewcommand{\k}{\lie{k}}

\newcommand{\g}{\lie{g}}

\newcommand{\on}{\operatorname}

\newcommand{\ad}{ \on{ad} }

\newcommand{\End}{ \on{End} }

\newcommand{\D}{ \mathcal{D} }

\newcommand\dirac{/\kern-1.2ex\partial} 
\newcommand\qu{/\kern-.7ex/} 

\renewcommand{\d}{{\mbox{d}}}
\newcommand{\ol}{\overline}

\newcommand\eps{\epsilon}
\newcommand\Om{\Omega}

\newcommand{\f}{\frac}

\newcommand{\p}{\partial}

\newcommand{\hh}{{\textstyle \f{1}{2}}}
\newcommand{\ti}{\tilde}

\newcommand\beqn{\begin{equation}}      
\newcommand\eeqn{\end{equation}}      
\newcommand{\ca}{\mathcal}
\newcommand{\wh}{\widehat}

\newcommand{\mf}{\mathfrak}
\newcommand{\beq}{\begin{eqnarray*}}
\newcommand{\eeq}{\end{eqnarray*}}
\newcommand{\Duf}{\on{Duf}}

\begin{document}

\title[]{On the Kashiwara-Vergne conjecture}

\author{A. Alekseev}
\address{University of Geneva, Section of Mathematics,
2-4 rue du Li\`evre, c.p. 64, 1211 Gen\`eve 4, Switzerland}
\email{alekseev@math.unige.ch}

\author{E. Meinrenken}
\address{University of Toronto, Department of Mathematics,
40 St George Street, Toronto, Ontario M5S 2E4, Canada }
\email{mein@math.toronto.edu}

\date{\today}

\begin{abstract}
Let $G$ be a connected Lie group, with Lie algebra $\g$. In 1977, Duflo
constructed a homomorphism of $\g$-modules $\on{Duf}\colon S(\g)\to
U(\g)$, which restricts to an algebra isomorphism on invariants.
Kashiwara and Vergne (1978) proposed a conjecture on the
Campbell-Hausdorff series, which (among other things) extends the
Duflo theorem to germs of bi-invariant distributions on the Lie
group $G$. 

The main results of the present paper are as follows. (1) Using a
recent result of Torossian (2002), we establish the Kashiwara-Vergne
conjecture for any Lie group $G$. (2) We give a reformulation of the
Kashiwara-Vergne property in terms of Lie algebra cohomology.  As a
direct corollary, one obtains the algebra isomorphism $H(\g,S(\g))\to
H(\g,U(\g))$, as well as a more general statement for
distributions. 
\end{abstract}
\subjclass{} 

\maketitle
\tableofcontents
\section{Introduction}
Let $G$ be a connected Lie group, with Lie algebra $\g$. Let $U(\g)$
denote the universal enveloping algebra, and $S(\g)$ the symmetric
algebra. The Duflo map is the isomorphism of $\g$-modules,
\begin{equation}\label{eq:duflomap}
 \on{Duf}\colon S(\g)\to U(\g),
\end{equation}
obtained by precomposing the symmetrization map with an infinite-order
differential operator $\wh{J^{1/2}}$ on the symmetric algebra (viewed
as the space of polynomials on $\g^*$).  Here $J\in C^\infty(\g)$ is
the analytic function
\begin{equation}\label{eq:J}
J(x)={\det}_\g\Big(\f{1-e^{-\ad_x}}{\ad_x}\Big),
\end{equation}
and $\wh{J^{1/2}}\colon S(\g)\to S(\g)$ is obtained from the Taylor
series expansion of $J^{1/2}$ by replacing the variables $x$ with
derivatives $\f{\p}{\p \mu}$. Duflo's theorem \cite{du:op} states that
the map $\on{Duf}$ restricts to an \emph{algebra} isomorphism on
invariants, $S(\g)^\g\to U(\g)^\g$. In \cite{ka:ca} Kashiwara and
Vergne proposed a conjecture regarding the Campbell-Hausdorff series,
which among other things generalizes the Duflo theorem to germs of
invariant distributions on the Lie group $G$. Here $S(\g)\subset
\D'_{\on{comp}}(\g)$ is identified with the subalgebra (under
convolution) of distributions on $\g$ supported at $0$,
while $U(\g)$ is identified with the subalgebra (under convolution) 
of distributions on $G$ supported at $e$. This generalization of 
Duflo's theorem was proved later in a series of papers
\cite{an:de,an:ko,an:co} using the diagrammatic technique of
Kontsevich \cite{ko:de}.  The KV conjecture itself was established for
$\g$ solvable \cite{ka:ca}, for $\g=\mf{sl}(2,\R)$ in \cite{ro:de}, and
in the case of $\g$ quadratic (that is, $\g$ carries an invariant
nondegenerate symmetric bilinear form) in \cite{ve:ce} (see also
\cite{al:ka}). 

There are several equivalent formulations of the KV
conjecture , one of which is as follows (see Section \ref{sec:review}
for details).
Let $\g_t$ denote the family of Lie algebras, obtained from $\g_1=\g$
by rescaling the Lie bracket as $[\cdot,\cdot]_t=t[\cdot,\cdot]$, and
let $m_t$ denote the family of products on $S(\g)$, induced from the
products on $U(\g_t)$ by the Duflo isomorphisms $\on{Duf}_t\colon
S(\g)\to U(\g_t)$. Thus, $m_t$ interpolates between the usual
commutative product $m_0$ on $S(\g)$ and the non-commutative product
coming from $U(\g)$. The product maps $m_t$ have a natural extension
to the space of distributions on $\g$, with compact support in a
sufficiently small open neighborhood $O$ of $0$,
\begin{equation}\label{eq:mt}
m_t: \D'_{\on{comp}}(O^2) \to \D'_{\on{comp}}(\g).
\end{equation}
Let $e_i$ be a basis of $\g^2=\g \oplus \g$, and let $L(e_i)$ be the
Lie derivatives relative to the adjoint action of $\g^2$ on
$\D'_{\on{comp}}(O^2)$. The Kashiwara-Vergne conjecture asserts the
existence of an analytic map $\beta=\sum \beta^i e_i: O^2 \to \g^2$,
vanishing at the origin, such that
\begin{equation} \label{eq:1}
\frac{dm_t}{dt}=- m_t \circ \sum_i \beta^i_t L(e_i),
\end{equation}
where $\beta^i_t(x,y)=t^{-1} \beta^i(tx,ty)$.  Equation \eqref{eq:1}
implies that the family of products $m_t$ is independent of $t$ when
restricted to germs of invariant distributions. 

In a recent paper, Torossian \cite{to:su} came very close to a proof
of the KV conjecture, again using Kontsevich diagram techniques. In
more detail, Torossian proved an analogue of Equation \eqref{eq:1},
for a certain family of maps $\ti{m}_t: \D'_{\on{comp}}(O^2) \to
\D'_{\on{comp}}(\g)$ with $\ti{m}_0=m_0, \ti{m}_1=m_1$. However,
$\ti{m}_t\not=m_t$, and indeed the products defined by $\ti{m}_t$ are
not associative.

Our first main result shows that one obtains a solution $\beta_t$ of
the KV equation \eqref{eq:1} from a solution $\ti{\beta}_t$ of
Torossian's modified KV property, by solving an interesting `zero
curvature' equation. This establishes the Kashiwara-Vergne conjecture
in full generality.

Our second result relates the KV property to Lie algebra
cohomology, as follows. For any $\g$-module $X$, let $H^\bullet(\g,X)$ denote the
Lie algebra cohomology of $\g$ with coefficients in $X$. By
functoriality, the Duflo map \eqref{eq:duflomap} induces a
homomorphism of graded vector spaces,
\begin{equation}\label{eq:hduflo}
H^\bullet(\g,S(\g))\to H^\bullet(\g,U(\g)).
\end{equation}
In \cite{sho:ts}, Shoikhet announced a joint result with Kontsevich
that the linear map \eqref{eq:hduflo} is an isomorphism of graded
algebras.  Shortly later, a proof was published by Pevzner-Torossian
\cite{pev:is}. In his paper, Shoikhet raises the question
\cite[Remark 3.2]{sho:ts} whether this result would follow from the KV
conjecture.  We show that this is indeed the case:  Let 
\begin{equation}M_t\colon 
\D'_{\on{comp}}(O^2)\otimes \wedge(\g^2)^*\to 
\D'_{\on{comp}}(\g)\otimes\wedge\g^*\end{equation}
denote the family of products on Chevalley-Eilenberg complexes,
induced by the maps $m_t$ \eqref{eq:mt}.  We prove that $\beta: O^2
\to \g^2$ solves the KV conjecture if and only if
\begin{equation}\label{eq:homform}
\f{d M_t}{d t}=-\big[\d,M_t\circ \sum_i \beta_t^i \otimes \iota(e_i) \big] .
\end{equation}
Here $\d$ is the Chevalley-Eilenberg differential, and $\iota(e_i)$
are contraction operators.  Equation \eqref{eq:homform} implies that
the product maps $M_t$ are chain homotopic. Hence, the induced map 
in cohomology is independent of $t$. 

The organization of this paper is as follows. In Section
\ref{sec:review}, we give three different formulations I,II,III of the
KV conjecture. The algebraic formulation III is the main version of
the conjecture, as stated in \cite[page 250]{ka:ca}.  The other two
versions are more geometric. They are implicit in \cite[Section
3]{ka:ca}, although we add some details to show that they are indeed
equivalent to III.  In Section \ref{sec:tor} we show that Torossian's
theorem \cite{to:su} implies the KV conjecture.  Section
\ref{sec:liecoh} gives yet another reformulation IV of the
KV conjecture, as a homotopy formula for
Chevalley-Eilenberg complexes. In this context, we formulate
a new conjecture extending the Kashiwara-Vergne approach
to higher homotopies.
 \\

\noindent{\bf Notational conventions.}
For any manifold $M$, we denote by $\mf{X}(M)$ the Lie algebra of 
vector fields (viewed as derivations of $C^\infty(M)$), and by 
$\Om^\bullet(M)$ the algebra of differential forms. 
We follow the convention that the flow $F_t$ of a 
time dependent vector field $X_t$ is given by 
$(X_t f)(F_t^{-1}(x))=\f{\p}{\p t}f(F_t^{-1}(x))$. The Lie derivative 
on $\Om(M)$ is thus given by $F_t^*\circ L_{X_t}=-\f{\p}{\p t} F_t^*$. 
A family of forms $\alpha_t$ is intertwined by the flow if 
$F_t^*\alpha_t$ is independent of $t$; in terms of the Lie derivative 
this means $(\f{\p}{\p t}-L(X_t))\alpha_t=0$. \\

\noindent{\bf Acknowledgments.} 
We would like to thank C. Torossian for explaining his work
\cite{to:su}. We are grateful to M. Vergne for comments, and to
H. Maire and E. Petracci for helpful discussions.  We are indebted to
the referees of this paper for valuable remarks. The research of
A.A. was supported in part by the Swiss National Science Foundation.

\section{Review of the Kashiwara-Vergne conjecture}\label{sec:review}
\subsection{Preliminaries}
Throughout this paper, $r_t^V$ will denote scalar multiplication by a
real parameter $t$ on some vector space $V$. We omit the
superscript if $V$ is clear from the context. We will frequently use
the following elementary fact. Suppose $R_t\colon W\to W,\ 0<t\le 1$
are automorphisms of a topological vector space, with $R_{t_1
t_2}=R_{t_1} \circ R_{t_2}$.  Assume that $w\in W$ be a smooth vector,
i.e. $w_t=R_t(w)$ is differentiable in $t$. Then the derivative
$\dot{w}_t=\f{\p w_t}{\p t}$ satisfies
\begin{equation}
 \dot{w}_t=t^{-1}R_t(\dot{w}_1).
\end{equation}
For example, suppose $O$ is a star-shaped open subset of a
vector space $V$, and take $W$ to be the space of $k$-forms on $O$. 
Let $\alpha\in\Om^k(O)$ and $a\in \R$. Then the derivative of 
$\alpha_t:=t^a r_t^* \alpha$ scales according to 
$\dot{\alpha}_t=t^{a-1}r_t^*\dot{\alpha}_1$. 

\subsection{The Duflo map on distributions}
Let $\g$ be a finite-dimensional Lie algebra over $\R$, and $G$ the corresponding
connected, simply connected Lie group. Let $\exp\colon \g\to G$ denote
the exponential map.  The Jacobian $J\in C^\infty(\g)$ of the
exponential map (using left trivialization of the tangent
bundle of $G$) is given by formula \eqref{eq:J}.  Let
$\D'_{\on{comp}}(\g)$ denote the commutative algebra of compactly
supported distributions on $\g$, with product given by convolution.
The symmetric algebra $S(\g)$ is identified with the subalgebra of
distributions supported at $0$. Similarly, let $\D'_{\on{comp}}(G)$
denote the non-commutative convolutions algebra of distributions on
$G$, and identify $U(\g)$ with distributions supported at the group
unit $e\in G$. Assume temporarily that $J$ admits a global square root
$J^{1/2}$, equal to $1$ at the origin. (This is not the case for all
$G$, but holds true, for instance, if $\g$ admits an invariant
non-degenerate symmetric bilinear form.) Then the Duflo map
$\Duf\colon S\g\to U\g$ extends to a homomorphism of $\g$-modules,
\begin{equation}\label{eq:duflo00}
 \Duf=\exp_*\circ J^{1/2}\colon \D'_{\on{comp}}(\g)\to
 \D'_{\on{comp}}(G),
\end{equation}
and again this is an algebra homomorphism on invariants. For a general
non-compact Lie group, this extension may not be very interesting
since there may be very few invariant distributions of compact
support. However, one obtains an interesting generalization by working
instead with \emph{germs} of invariant distributions. Let 
$\mf{D}'(\g,A)$ denote the space of germs (at the origin $0\in\g$)
of distributions on $\g$, with asymptotic support \footnote{Let $x$ be
  a manifold, and $\mf{u}$ the germ of a distribution at $x\in M$.
  The {\em (asymptotic) support} of $\mf{u}$ is the closed subset
  $\on{supp}_x(\mf{u})\subset T_xM$, defined as follows: Let
  $\gamma(t)$ denote a smooth curve with $\gamma(0)=x$. Then the
  tangent vector $\dot{\gamma}(0)\not\in \on{supp}_x(\mf{u})$ if and
  only if $\gamma(t)\not\in \on{supp}(\mf{u})$ for $t>0$ sufficiently
  small.}  in the closed cone $A\subset \g$. Let $\mf{D}'(G,A)$ be the
space of germs of distributions at $e\in G$, given as the image of
$\D'(\g,A)$ under $\exp_*$. For closed cones $A_1,A_2\subset\g$ with
$A_1\cap -A_2=\{0\}$, there are well-defined convolution products
\begin{equation}\label{eq:alge}
\mf{D}'(\g,A_1)\otimes  \mf{D}'(\g,A_2)\to \mf{D}'(\g,A_1+A_2),
\end{equation}
and 
\begin{equation}\label{eq:group}
\mf{D}'(G,A_1)\otimes  \mf{D}'(G,A_2)\to \mf{D}'(G,A_1+A_2).
\end{equation}
For instance, these products are well-defined if $A_1=\{0\}$ or
$A_2=\{0\}$.  The natural extension of Duflo's theorem asserts that
the Duflo map intertwines \eqref{eq:alge} (restricted to
$\g$-invariants), with \eqref{eq:group} (restricted to
$\g$-invariants).  This extension of Duflo's theorem was put forward
by Kashiwara-Vergne \cite{ka:ca} in their paper, and proved for the
case of solvable Lie algebras. The general case was proved many
years later in \cite{an:de,an:ko,an:co}. 

\subsection{The Kashiwara-Vergne approach}
The Kashiwara-Vergne approach is based on the idea that the Duflo
theorem (and its extension to germs of distributions) should follow
from a more general property of the Campbell-Hausdorff series
\begin{equation} 
\Phi(x,y):=\log(\exp(x)\exp(y))=x+y+\hh [x,y]+\cdots\ .\end{equation}
Recall that the function $\Phi(x,y)$ is well-defined and analytic for
$x,y$ in a sufficiently small neighborhood of the origin. To be
specific, let $U\subset \g$ be a star-shaped neighborhood of the
origin such that the exponential map is a diffeomorphism over $U$, 
and let $O\subset U$ be a
smaller star-shaped open neighborhood such that $\exp(x)\exp(y)\in
\exp(U)$ for all $x,y\in O$. Then $\Phi\colon O^2\to\g$ is
well-defined and takes values in $U$.
Since $J>0$ on $U$, 
the square root $J^{1/2}$ is a well-defined analytic function on $U$.
Under the Duflo map, the convolution product on $\D'_{\on{comp}}(G)$ gives rise to a
non-commutative product on distributions on $\g$ with support in $O$. 
In terms of the map $\Phi$ this product can be written, 
\begin{equation}
 m=\Phi_* \circ \kappa\colon \D'_{\on{comp}}(O^2)\to \D'_{\on{comp}}(\g)
\end{equation}
where $\kappa\in C^\infty(O^2)$ is the function, 
\begin{equation}
\kappa(x,y)=\f{J^{1/2}(x) J^{1/2}(y)}{J^{1/2}(\Phi(x,y))}.
\end{equation}
Denote by $\g_t$ the Lie algebra, obtained from $\g$ by rescaling the
Lie bracket, 
\begin{equation} [\cdot,\cdot]_t=t[\cdot,\cdot].\end{equation}
Let $G_t$ denote the connected, simply connected Lie group with Lie
algebra $\g_t$.  For $t\not=0$, the map $r_t\colon\g\to \g$ defines a
Lie algebra isomorphism $\g_t\to \g_1$, while $\g_0$ is the vector
space $\g$ with the zero bracket.  The exponential map $\exp_t\colon
\g_t\to G_t$ is a diffeomorphism over $U_t=r_{t^{-1}} U$, and the
counterpart $\Phi_t\colon O_t\to \g$ of the map $\Phi$ is well-defined
over the rescaled neighborhood 
$O_t=r_{t^{-1}} O$. One has $\Phi_t=t^{-1}r_t^*\Phi$ for
$t\not=0$, i.e.
\begin{equation} \Phi_t(x,y)=x+y+\textstyle{\f{t}{2}}[x,y]+\cdots,\end{equation}
while $\Phi_0(x,y)=x+y$.  The Jacobian of $\exp_t$ is $J_t=r_t^*J$.
From the family of Lie brackets $[\cdot,\cdot]_t$, we obtain a family
of multiplication maps $m_t=(\Phi_t)_*\circ \kappa_t$ with $
\kappa_t=r_t^*\kappa$. For $t\not=0$,
\begin{equation} m_t=(r_{t^{-1}})_*\circ m_1\circ (r_t)_*.\end{equation}
The map $m_0$ is the commutative product given by convolution on the
vector space $\g$. The $t$-derivative $\dot{m}_t=\f{\p
m_t}{\p t}$ scales according to 
\begin{equation}
 \dot{m}_t=t^{-1}(r_{t^{-1}})_*\circ \dot{m}_1\circ (r_t)_*.\end{equation}
The generalized Duflo theorem is equivalent to the statement that the derivative 
$\dot{m}_t$ vanishes on $\D'(\g\times \g,A_1\times A_2)^{\g\times\g}$, where
$A_1\cap(-A_2)=\{0\}$. Kashiwara-Vergne formulated the following
stronger property of $m_t$. 

\subsection{The Kashiwara-Vergne conjecture}
Let $e_i$
be a basis of $\g^2$, and let $L_i=L(e_i)$ be the Lie derivatives for
the
$\g^2$-action on $\D'_{\on{comp}}(O^2)$.  \\

\noindent {\bf Kashiwara-Vergne conjecture, version I}: 
{\em Taking $O$ smaller if necessary, there exists a $\g$-equivariant 
smooth function $\beta=\sum_i \beta^i e_i\colon O^2\to \g^2$, with $\beta(0,0)=0$, such that
  $\beta_t=t^{-1}r_t^*\beta$ satisfies,
\begin{equation}\label{eq:KV1}
\f{\p m_t}{\p t}= -m_t \circ \sum_i \beta_t^i L_i.
\end{equation}
}\\
Duflo's theorem, as well as 
its extension to germs of invariant distributions, are a direct consequence of 
\eqref{eq:KV1} since the right hand side vanishes on invariants.

\begin{remark}
It is not difficult to show that $\dot{m}_t$ can be written in the
form
\begin{equation}\label{eq:pt}
 \dot{m}_t=m_t\circ P_t
\end{equation}
where $P_t$ are first order differential operators on
$\D'_{\on{comp}}(O^2)$. In light of Duflo's theorem, the idea that one
can take $P_t$ of the form $\sum_i \beta_t^i L_i$ may therefore appear
rather natural. Note also that the operator $\sum_i \beta_t^i L_i$
on $\D'(O^2)$ is different, in general, from the Lie derivative $L(X^{\beta_t})$ 
in the direction of the vector field $X^{\beta_t}$ generated by 
$\beta_t$, since on distributions $L(X^{\beta_t})=\sum_i L_i
\circ \beta_t^i$.
\end{remark}
 
%
%
%
%
%
Another version of the conjecture states the KV property in terms of the 
derivatives of the function $\Phi_t$ and $\kappa_t$.\\ 

\noindent {\bf Kashiwara-Vergne conjecture, version II}: 
{\em Taking $O$ smaller if necessary, there exists a $\g$-equivariant 
smooth function
  $\beta\colon O^2\to \g^2$, with $\beta(0,0)=0$, such that
  $\beta_t=t^{-1}r_t^*\beta$ satisfies the following two equations:
\begin{equation} \label{eq:KV2a}
\f{\p \Phi_t}{\p t}-L(X^{\beta_t}) \Phi_t=0,\ \ 
\end{equation}
\begin{equation}\label{eq:KV2b}
\f{\p \kappa_t}{\p t}-\big(\sum_i L_i \circ \beta^i_t \big) \kappa_t=0.
\end{equation}
%
}\\

\begin{remark}
The geometric interpretation of the first equation is that the vector
field $X^{\beta_t}$ intertwines the product maps, $\Phi_t$. To
interpret the second equation, choose a translation invariant volume
form $\Gamma_0$ on $\g^2$, and consider the family of volume forms
$\Gamma_t=\kappa_t \Gamma_0$. Then the second equation is equivalent to
\[ \begin{split}
\dot{\Gamma}_t&=L(X^{\beta_t})\Gamma_t-\sum_i \beta_t^i \kappa_t L_i \Gamma_0\\
&=L(X^{\beta_t})\Gamma_t -\sum_i \beta_t^i \on{tr}_{\g\times\g}(e_i) \Gamma_t
\end{split}
\]
where we have used that $L(X^{\beta_t})=\sum_i L_i\circ
\beta_t^i$ on volume forms. If $\g$ is \emph{unimodular} (i.e.,
$\on{tr}_\g(\ad(x))=0$ for all $x$), the equation says that the vector
field $X^{\beta_t}$ intertwines the volume forms $\Gamma_t$, while 
in the general case $\Gamma_t$ is a relative invariant
for the unimodular character.
\end{remark}

The third version of the Kashiwara-Vergne conjecture is more
algebraic, and no longer involves the variable $t$. 
This is the `main' version of the conjecture, as formulated in  \cite{ka:ca}.
\\

\noindent {\bf Kashiwara-Vergne conjecture, version III}: 
{\em Taking $O$ smaller if necessary, 
there exist smooth $\g$-equivariant functions 
$\beta^1,\beta^2\colon O^2\to \g$, vanishing at
$(0,0)$, such that the following two equations are satisfied:
\begin{equation}\label{eq:KV3a}
(1-e^{-\ad_x})\beta^1(x,y)+(e^{\ad_y}-1)\beta^2(x,y)=x+y-\log(\exp(y)\exp(x)),
\end{equation}
\begin{equation}\label{eq:KV3b}
-\on{tr}\Big(\ad_x \circ \partial_x \beta^1(x,y)+\ad_y \circ \partial_y \beta^2(x,y)\Big)
=\hh \on{tr}\Big(\f{\ad_x}{e^{\ad_x}-1}+\f{\ad_y}{e^{\ad_y}-1}
-\f{\ad_z}{e^{\ad_z}-1}-1\Big),
\end{equation}
with $z=\log(\exp(x)\exp(y))$. Here $\on{tr}$ denotes the trace in the
adjoint representation,  
$\partial_x \beta^1(x,y)\in \on{End}(\g)$ is the linear map 
\begin{equation}\xi\mapsto\f{\p}{\p u}\Big|_{u=0}\beta^1(x+u\xi,y),\end{equation}
and similarly for $\partial_y \beta^2(x,y)$.
}\\

\begin{remark}\label{rem:universal}
  A \emph{Lie polynomial} in variables $x_1,\ldots,x_l$ 
  is an element of the free Lie algebra $\mf{l}=\mf{l}(x_1,\ldots,x_l)$ 
  with generators $x_i$. Recall that the free Lie
  algebra has a grading $\mf{l}=\bigoplus_k \mf{l}^k$, where the
  degree $k$ component consists of linear combinations of \emph{Lie monomials}
  $[z_{i_1},[z_{i_2},\ldots,[z_{i_{k-1}},z_{i_k}],\ldots ]]$ where
  $z_j\in \{x_1,\ldots,x_l\}$. A \emph{Lie series} is an element of the
  degree completion $\wh{\mf{l}}=\prod_k \mf{l}^k$. If $x_i$ are elements
  of a given Lie algebra $\g$, one obtains a Lie algebra homomorphism
  $\mf{l}\to \g$ taking $x_i\in\mf{l}$ to $x_i\in \g$.  
  
  The Campbell-Hausdorff theorem (see e.g.  \cite[Chapter I.3]{go:fo})
  states that if $\dim\g<\infty$, the function
  $\Phi(x,y)=\log(\exp(x)\exp(y))$ is a Lie series in $x,y$. Dynkin's
  formula gives this series explicitly as a \emph{universal} element
  in $\wh{\mf{l}}(x,y)$ (i.e. not depending on the Lie algebra).  It
  is hence natural to require that similarly, $\beta^1,\beta^2$ are
  analytic functions in $x,y$ given as universal Lie series. 
 \end{remark}
\subsection{Equivalence of versions I-III of the KV conjecture}
The implication $III\Rightarrow II\Rightarrow I$ was shown in
\cite[Section 3]{ka:ca}. Closer examination of their argument shows
that the three statements are in fact equivalent: 
\subsubsection{Equivalence $I\Leftrightarrow II$}
We have to show that $\dot{m}_t= m_t \circ \sum_i \beta_t^i
L_i$ if and only if $\beta_t$ solves the pair of equations
\eqref{eq:KV2a},\eqref{eq:KV2b}. In fact, this equivalence does not
involve the particular form of $\Phi_t,\kappa_t$. Thus let 
$O\subset \g$ be an open neighborhood of the origin, and let $\Phi_t:
O^2 \to \g, \kappa_t: O^2 \to \R$ be \emph{arbitrary} families of
smooth functions. Define a family of maps $\D'_{\on{comp}}(O^2) \to
\D'_{\on{comp}}(\g)$ by $m_t=(\Phi_t)_* \circ \kappa_t$. Let $\beta_t:
O^2 \to \g^2$ be a smooth family of functions.
\begin{lemma}\label{lem:tran}
The equation 
$\dot{m}_t =- m_t \circ \sum_i \beta_t^i L_i$ holds if and only if 
\begin{equation} \label{eq:phikappa}
\dot{\Phi}_t - L(X^{\beta_t}) \Phi_t=0,\ \ 
\dot{\kappa}_t - \big(\sum_i L_i \circ \beta^i_t \big) 
\kappa_t=0.
\end{equation}
\end{lemma}
\begin{proof}
Recall that linear
combinations of delta-distributions are dense in the space of
distributions. Hence $\dot{m}_t = -m_t \circ \sum_i
\beta^i L_i$ holds if and only if it holds on
$\delta_{(q,p)}$, for all $(q,p)\in O^2$. Let us work out the resulting
identities. Taking the $t$-derivative of 
\begin{equation} m_t(\delta_{(q,p)})=\kappa_t(q,p) \delta_{\Phi_t(q,p)}, \end{equation}
we obtain
\begin{equation}\label{eq:tder}
 \dot{m}_t(\delta_{(q,p)})=\dot{\kappa}_t(q,p)\delta_{\Phi_t(q,p)}
-\kappa_t(q,p) \sum_a \dot{\Phi}_t^a(q,p) 
\f{\p }{\p x^a}\delta_{\Phi_t(q,p)}.\end{equation}
Here $e_a$ is a basis of $\g$ and $x^a$ are the corresponding
coordinates. On the other hand, 
\beq 
\big(m_t \circ \sum_i \beta_t^i L_i\big) (\delta_{(q,p)})&=& 
(\Phi_t)_*\Big(\sum_i {\kappa}_t \beta_t^i L_i\,
\delta_{(q,p)}\Big)\\
&=& (\Phi_t)_* \Big(-\sum_i \big(L_i({\kappa}_t\beta_t^i)\big) \delta_{(q,p)}
+\sum_i ({\kappa}_t\beta_t^i)(q,p) L_i \, \delta_{(q,p)}\Big)\\
&=&-\sum_i L_i(\kappa_t\beta_t^i )(q,p) \delta_{\Phi_t(q,p)} +\sum_i
(\kappa_t\beta_t^i )(q,p) (\Phi_t)_* (L_i\,\delta_{(q,p)})\\ &=&-\sum_i
L_i(\kappa_t\beta_t^i)(q,p) \delta_{\Phi_t(q,p)} +\sum_{ia}
(\kappa_t\beta_t^i\, L_i(\Phi_t^a))(q,p) \f{\p }{\p
x^a}\delta_{\Phi_t(q,p)}.  \eeq
Comparing the coefficients of $\delta_{\Phi(p,q)},\ \f{\p }{\p
x^a}\delta_{\Phi_t(q,p)}$ with those in Equation \eqref{eq:tder}, we
find that $\dot{m}_t =- m_t \circ \sum_i \beta^i L_i$ is equivalent to 
\begin{equation} \dot{\Phi}_t^a=\sum_i \beta_t^i L_i(\Phi_t^a),\ \  
\dot{\kappa}_t=\sum_i L_i(\kappa_t\beta_t^i ) ,
\end{equation}
which is \eqref{eq:phikappa}. 
\end{proof}

\subsubsection{Equivalence $II\Leftrightarrow III$.}
Let $A(x,y), A'(x,y)$ denote the left-, right-hand side of \eqref{eq:KV3a}, and
write $A_t=t^{-1}r_t^*A$ and
$A'_t=t^{-1}r_t^* A'$. In the proof of \cite[Lemma 2.3]{ka:ca},
Kashiwara-Vergne compute 
\[ \begin{split}
\dot{\Phi}_t(x,y)&=
\f{1}{t}\f{\p}{\p \eps}\Big|_{\eps=0}\log\Big(e^{tx}e^{\eps A_t'(x,y)}e^{ty}\Big)\\
\\ 
-(L(X^{\beta_t})\Phi_t)(x,y)&=\f{1}{t}\f{\p}{\p
  \eps}\Big|_{\eps=0}\log\Big(e^{tx}e^{\eps A_t(x,y)}e^{ty}\Big).
\end{split}\]
It follows that the first condition \eqref{eq:KV2a} of version II is equivalent to
$A_t=A'_t$, hence to $A=A'$ which is the first condition
\eqref{eq:KV3a} of version III. 

Assuming \eqref{eq:KV2a}, we now show that \eqref{eq:KV2b} 
is equivalent to \eqref{eq:KV3b}.  
Let $B(x,y)$,
$B'(x,y)$ denote the left-, right-hand side of Equation
\eqref{eq:KV3b}, and let $B_t=t^{-1}r_t^*B$, $B'_t=t^{-1}r_t^*B'$.
Then 
\begin{equation}\label{eq:BTT}
B_t=-\sum_i (L_i\beta_t^i),
\end{equation}
and Equation \eqref{eq:KV2b} may be written in the form
\begin{equation}
0=\dot{\kappa}_t-L(X^{\beta_t})\kappa_t+\sum_i (L_i\beta^i_t)\kappa_t=
\big(\f{\p}{\p t}-L(X^{\beta_t})\big)\kappa_t-B_t \kappa_t.
\end{equation}
From \eqref{eq:KV2a}, we have
\[ \big(\f{\p}{\p t}-L(X^{\beta_t})\big)J_t(\Phi_t(x,y))=\dot{J}_t(\Phi_t(x,y)).\]
On the other hand,
$L(X^{\beta_t})J_t(x)=L(X^{\beta_t})J_t(y)=0$ since $J$ is
$\ad$-invariant. 
Recalling the formula \cite[Lemma 3.3]{ka:ca}
\begin{equation} J_t(z)^{-1}\dot{J}_t(z)=\on{tr}_\g\Big(\f{t\ad_{z}}{e^{t\ad_z}-1}-1\Big),\ \ z\in\g\end{equation}
this gives, 
\begin{equation}\dot{\kappa}_t-L(X^{\beta_t})\kappa_t
=B'_t \kappa_t.\end{equation}
We thus obtain $B_t=B_t'$, i.e. $B=B'$.

\section{Torossian's theorem implies the KV conjecture}
\label{sec:tor}
\subsection{Torossian's theorem}
In his paper, Torossian \cite[page 601]{to:su} proves the following
statement, similar to version II of the KV conjecture.
\begin{theorem}[Torossian]\label{th:tor}
There are families of smooth functions 
\begin{equation}
 \ti{\Phi}_s\colon O^2\to \g,\ \ti{\kappa}_s\colon O^2\to\R,\  \
\ti{\beta}_s\colon O^2\to \g^2,\ \ \ s\in[0,1]\end{equation}
with the following properties:
\begin{enumerate}
\item
The Taylor expansion of $\ti{\Phi}_s(x,y)$ at the origin $(0,0)$ is 
of the form $\ti{\Phi}_s(x,y)=x+y+\cdots$.  Furthermore, 
$\ti{\kappa}_s(0,0)=1,\ \ \ti{\beta}_s(0,0)=0$, 
\item
For $s=0$, 
\[
  \ti{\Phi}_{0}(x,y)=x+y,\ \ \ti{\kappa}_0(x,y)=1,\ \ \ 
\]
while for $s=1$, 
\[
\ti{\Phi}_1(x,y)=\Phi(x,y),\ \ \ti{\kappa}_1(x,y)=\kappa(x,y).
\]
\item
The following two equations are satisfied:  
\[
\dot{\ti{\Phi}}_s - L(X^{\ti{\beta}_s}) \ti{\Phi}_s=0, \ \ \ 
\dot{\ti{\kappa}}_s - \big(\sum_i L_i \circ 
\ti{\beta}^i_s \big) 
\ti{\kappa}_s=0.
\]
\end{enumerate}
\end{theorem}
The functions $\ti{\Phi}_s,\ \ti{\kappa}_s,\ \ti{\beta}_s$ are
constructed using the diagrammatic technique of Kontsevich
\cite{ko:de}. See \cite{to:su} for details.

Using $\ti{\Phi}_s$ and $\ti{\kappa}_s$ one can define maps 
$\ti{m}_s=(\ti{\Phi}_s)_*\circ \ti{\kappa}_s$ as before.  However,
$\ti{m}_s\not=m_s$, and in fact the `products' $\ti{m}_s$ are not
associative, in general. Nonetheless, according Lemma \ref{lem:tran}
one still has 
\begin{equation}
\f{\p \ti{m}_s}{\p s}=-\ti{m}_s\circ \sum_i \ti{\beta}_s^i L_i,
\end{equation}
which as before implies the Duflo theorem and its generalization to
germs of invariant distributions.

\subsection{Torossian $\Rightarrow$ KV}
We will now show how to obtain a solution $\beta_t$ of the original KV problem
from Torossian's function $\ti{\beta}_s$. Let us first state the
relevant argument 
in a more general
setting. Let $\k$ be a Lie algebra, acting on a manifold
$M$. The vector fields defining this action are a Lie algebra
homomorphism $\k\to \mf{X}(M),\ \xi\mapsto X^\xi$. For $\beta\in \k_M:=C^\infty(M,\k)$ define
$X^\beta\in\mf{X}(M)$ by $X^\beta(x)=X^{\beta(x)}(x)$. Then
\begin{equation}
\k_M\to \mf{X}(M),\ \beta\mapsto X^\beta
\end{equation}
is a Lie algebra homomorphism for the following Lie bracket on $\k_M$,
\begin{equation}\label{eq:newbracket}
[\beta,\gamma](x)=(L(X^\beta)\gamma)(x)-(L(X^\gamma)\beta)(x)+
[\beta(x),\gamma(x)].
\end{equation}
The following Lemma will be proved in the appendix.
\begin{lemma}\label{lem:zero}
Suppose $\gamma_{s,t}\in\k_M$ is a given smooth 2-parameter family,
with $\gamma_{s,t}(x_0)=0$ for a given point 
$x_0\in M$. Replacing $M$
by a smaller neighborhood of $x_0$ if necessary, there exists a 
unique 2-parameter family $\beta_{s,t}\in \k_M$ with
$\beta_{0,t}=0$, such that the following zero 
curvature equation is satisfied: 
\begin{equation} \label{eq:flat1}
\frac{\partial \beta_{s,t}}{\partial s} - 
\frac{\partial \gamma_{s,t}}{\partial t}
+ [\beta_{s,t}, \gamma_{s,t}]= 0. 
\end{equation}
\end{lemma}
Let us describe another general feature of the Lie algebra
$\k_M$. For $\beta\in\k_M$, consider the first order differential operator
on $\D'(M)$ (cf. \eqref{eq:KV2b})
\begin{equation} 
V(\beta)=\sum_i \beta^i\circ L(e_i) 
\end{equation}
The operator $V(\beta)$ is different from the Lie derivative
$L(X^\beta)$, which equals $\sum_i L_i\circ \beta^i$ on
distributions. However, one verifies that the difference
$L(X^\beta)-V(\beta)=\tau(\beta)$, where
\begin{equation} \tau\colon \k_M \to C^\infty(M),\ 
\beta\mapsto \sum_i (L_i \beta^i),
\end{equation}
is a Lie algebra cocycle: 
$\tau([\beta,\gamma])=L(X^\beta)\tau(\gamma)-L(X^\gamma)\tau(\beta)$.
This gives
\begin{equation} V([\beta,\gamma])=[V(\beta),V(\gamma)].\end{equation}
(Note that for the function $\beta$ from the Kashiwara-Vergne conjecture
(version III), the cocycle $\tau(\beta)$ appears as the left hand side of
\eqref{eq:KV3b} (see also \eqref{eq:BTT}).)

\begin{theorem}\label{th:suppl}
  Any solution $\ti{\beta}_s$ of Torossian's modified KV problem (cf.
  Theorem \ref{th:tor}) determines a solution $\beta$ of the KV
  problem.
\end{theorem}

\begin{proof}
Introduce a 2-parameter family of maps
$\ti{\beta}_{s,t}=r_t^*\ti{\beta}_s,\ \ s,t\in [0,1]$.
Since $\ti{\beta}_s(0,0)=0$, we have $\ti{\beta}_{s,t}(0,0)=0$. 
Similarly, rescale $\ti{\Phi}_{s,t}=t^{-1}r_t^*\ti{\Phi}_s$ and
$\ti{\kappa}_{s,t}=r_t^*\ti{\kappa}_s$. For each $s,t$ the `product' 
$\ti{m}_{s,t}=(\ti{\Phi}_{s,t})_*\circ \ti{\kappa}_{s,t}$ satisfies
the equation, 
$$
\f{\p \ti{m}_{s,t}}{\p s}=-\ti{m}_{s,t}\circ \sum_i \ti{\beta}_{s,t}^i L_i .
$$
Now let $\beta_{s,t}$ be the  2-parameter family of maps $\beta_{s,t}: O^2
\to \g^2,\ \ \beta_{0,t}=0$ obtained by solving the zero curvature
equation, Lemma \ref{lem:zero}, with $\gamma_{s,t}=\ti{\beta}_{s,t}$. 
The scaling property $\ti{\beta}_{s,t}=r_t^*\ti\beta_s$ implies the 
scaling property for the derivative, 
$\f{\p\ti{\beta}_{s,t}}{\p t}=t^{-1} 
r_t^*\big(\f{\p \ti{\beta}_{s,u}}{\p u}|_{u=1}\big)$, and hence (by Equation 
\eqref{eq:flat1}) 
$\beta_{s,t}= t^{-1} r_t^*\beta_{s,1}$.
We compute, 
\[
\begin{split}
\lefteqn{\f{\p}{\p s} \left( \f{\p \ti{m}_{s,t}}{\p t} + \ti{m}_{s,t} \circ 
V(\beta_{s,t}) \right)}\\
& =  \f{\p}{\p t}\f{\p\ti{m}_{s,t}}{\p s} +\f{\p \ti{m}_{s,t}}{\p s} 
\circ V(\beta_{s,t}) +
\ti{m}_{s,t} \circ V\Big( \f{\p \beta_{s,t}}{\p s}\Big) \\
& =  -\f{\p}{\p t} \Big( \ti{m}_{s,t} \circ V(\ti{\beta}_{s,t})\Big) 
-\ti{m}_{s,t} \circ V(\ti{\beta}_{s,t}) \circ V(\beta_{s,t}) 
+ \ti{m}_{s,t} \circ V\Big( \f{\p \beta_{s,t}}{\p s}\Big) \\
& = -\f{\p\ti{m}_{s,t}}{\p t} \circ V(\ti{\beta}_{s,t}) - \ti{m}_{s,t} \circ 
\left( V\Big( \f{\p \ti{\beta}_{s,t}}{\p t}\Big) - V\Big( 
\f{\p \beta_{s,t}}{\p s}\Big) + 
V(\ti{\beta}_{s,t}) \circ V(\beta_{s,t}) \right) \\
& =  -\left(\f{\p\ti{m}_{s,t}}{\p t}  +
\ti{m}_{s,t} \circ V(\beta_{s,t}) \right) \circ V(\ti{\beta}_{s,t})
\end{split}
\]
where in the last step we used \eqref{eq:flat1} with $\gamma_{s,t}=\ti{\beta}_{s,t}$. 
This equality is an ordinary first order differential equation (in
$s$) for the map
\begin{equation}\label{eq:mst1}
\ca{R}_{s,t}:=\f{\p\ti{m}_{s,t}}{\p t}  +\ti{m}_{s,t} \circ V(\beta_{s,t}).
\end{equation}
Since $\ti{m}_{0,t}=m_0$ and $\beta_{0,t}=0$, one has the initial condition 
$\ca{R}_{0,t}=0$.  Hence, by uniqueness of solutions of ordinary differential 
equations it follows that $\ca{R}_{s,t}=0$ for all $s$. 
In particular, $\ca{R}_{1,t}=0$ gives the desired  equation 
\[ \f{\p m_t}{\p t}+m_t\circ V(\beta_t)=0\]
where $\beta_t=\beta_{1,t}$ has the scaling property 
$\beta_t=t^{-1}r_t^*\beta$. 
\end{proof}

\subsection{Analyticity}
As pointed out in Remark \ref{rem:universal}, it is natural to impose the additional condition
on the KV solution $\beta$, that its components $\beta^1,\beta^2$ 
are given by  (universal) analytic Lie series. We will now verify that the
solution $\beta$ constructed in the last Section does indeed have 
this property. 

We introduce the following terminology. Let
$\mf{l}=\bigoplus_{n=1}^\infty \mf{l}^n$ be the free Lie algebra on
generators $x,y$, and $\hat{\mf{l}}$ its degree completion. 
An element $l\in \mf{l}^n$ will be called a \emph{Lie word of length $n$}
if it is obtained by bracketing $n$ elements
$x_{i_1},\ldots,x_{i_n}\in \{x,y\}$, in any order. For instance,
$[y,[[x,y],x]]$ is a Lie word of length $4$.
\begin{definition}
Given a Lie series $l=\sum_{n=1}^\infty l_n\in\hat{\mf{l}}$, let  $C_n(l)$ denote the
smallest possible  $\sum_i |c_{n,i}|$ for presentations 
\begin{equation}
l_n=\sum_{i} c_{n,i} l_{n,i}
\end{equation}
where $l_{n,i}$ are Lie words of length $n$.  The Lie series $l$
is called {\em strongly convergent} if $C_n(l)\le D^n$ for some $D>0$.
\end{definition}
The Campbell-Hausdorff series is an example of a strongly convergent
Lie series.  For any finite-dimensional Lie algebra $\g$, a strongly
convergent Lie series $l$ defines an analytic function $l(x,y)$ with
values in $\g$ on some neighborhood of the origin in $\g \times
\g$. Strongly convergent Lie series form a Lie subalgebra of
$\hat{\mf{l}}$.

Motivated by \eqref{eq:newbracket}, we introduce a new bracket on
$\hat{\mf{l}}\times\hat{\mf{l}}$, as follows. For any  $\beta=(\beta^1,\beta^2)$, let  
$L_\beta$ denote the derivation of $\hat{\mf{l}}$ given on generators by 
$L_{\beta}x=[\beta^1,x],\ L_\beta y=[\beta^2,y]$. If
$\gamma=(\gamma^1,\gamma^2)$ is a pair of Lie series, we set 
$L_\beta\gamma=(L_\beta \gamma^1,L_\beta\gamma^2)$. Now put 
\begin{equation}
[\beta,\gamma]=L_\beta \gamma-L_\gamma\beta+[\beta,\gamma]_0.
\end{equation}
where $[\cdot,\cdot]_0$ denotes the original (componentwise) bracket.
The Jacobi identity for the new bracket is a simple corollary
of the derivation property for $L_\beta$. 

Let us write $\ad_\beta=[\beta,\cdot]$, and
denote by $C_n(\beta)$ the maximum of $C_n(\beta^1)$ and
$C_n(\beta^2)$.  
\begin{lemma}
Suppose $\beta,\gamma\in \hat{\mf{l}}\times\hat{\mf{l}}$ with
$C_n(\beta),C_n(\gamma)\le D^n$ for some $D>0$. Then 
$C_n([\beta,\gamma])\le n^2 D^n$. More generally, If
$C_n(\beta_i),C_n(\gamma)\le D^n$ then 
\begin{equation}\label{eq:estimate}
 C_n(\ad_{\beta_k}\cdots \ad_{\beta_1} \gamma)
\le \f{(2 n^2)^k}{(2k-1)!!} \,D^n.
\end{equation}
\end{lemma}
Here $(2k-1)!!=(2k-1)(2k-3)\cdots 1$ for all $k>0$.  
\begin{proof}
Note that if $\beta,\gamma$ are homogeneous of degree $n_1,n_2$, then 
$[\beta,\gamma]$ is homogeneous of degree $n=n_1+n_2$, and 
\[  C_n([\beta,\gamma])\le (n+1) C_{n_1}(\beta)C_{n_2}(\gamma).\]
This follows since $C_n(L_\beta\gamma)\le n_2 C_{n_1}(\beta)C_{n_2}(\gamma)$, 
and similarly $C_n(L_\gamma\beta)\le n_1 C_{n_1}(\beta)C_{n_2}(\gamma)$, while  
$C_n([\beta,\gamma]_0)\le C_{n_1}(\beta)C_{n_2}(\gamma)$.
We now obtain \eqref{eq:estimate} by induction on $k$. The case $k=0$
is trivial, and for $k>0$ we find, using the induction hypothesis for 
$k-1$,
%
%
%
\[ \begin{split}
C_n(\ad_{\beta_k}\cdots \ad_{\beta_1} \gamma)
&\le (n+1) \sum_{j=1}^{n-1} C_{n-j}(\beta_k)C_{j}(\ad_{\beta_{k-1}}\cdots
\ad_{\beta_1} \gamma)\\
&\le (n+1) \f{1}{(2k-3)!!}\sum_{j=1}^{n-1}
(2 j^2)^{(k-1)}\ D^n \\
&\le 2^{k-1} (n+1) \f{n^{2k-1}}{(2k-1)!!}\ D^n  \le \f{2^k n^{2k}}{(2k-1)!!} \,D^n.
\end{split}\]
(In the last line, we used $\sum_{j=1}^{n-1} j^{2k-2}\le \int_{1}^n
j^{2k-2}\ d j$.)
\end{proof}

We will need one more remark. Suppose $l=\sum_{n=1}^\infty l_n\in\hat{\mf{l}}$ is a strongly
convergent Lie series. Then $l_t=r_t^*l=\sum_{n=1}^\infty t^{n}
l_n$ is strongly convergent for all $0\le t\le 1$, and since
$t^{n}\le 1$ one has the uniform estimate $C_n(l_t)\le
C_n(l)\le D^n$. Similarly, $\dot{l}_t=\f{\p l_t}{\p t}$ is strongly convergent, with 
a uniform  estimate on $C_n(\dot{l}_t)$. 

\begin{proposition}
Let $\gamma_s\in\hat{\mf{l}}\times\hat{\mf{l}}$, with coefficients depending
continuously on $s\in [0,1]$, and such that $C_n(\gamma_s)\le D^n$. Let
$\gamma_{s,t}=r_t^*\gamma_s$. Then the formal sum 
\begin{equation}\label{eq:sum}
 \beta_{s,t}=\sum_{k=0}^\infty 
\int_{0\leq s_0 \leq s_1 \dots \leq s_k \leq s} ds_0  \dots ds_k\ 
\big(\ad_{\gamma_{s_k,t}} \cdots \ad_{\gamma_{s_1,t}} \f{\p \gamma_{s_0,t} }{\p t}\big)
\end{equation}
defines a strongly convergent Lie series for $(s,t) \in [0,1]^2$, 
with $\beta_{0,t}=0$. 
Moreover, $\beta_{s,t}=t^{-1}r_t^*\beta_s$ and the following zero
curvature equation holds true,
\begin{equation} \label{eq:flat}
\partial_s \beta_{s,t} - \partial_t \gamma_{s,t} + [\beta_{s,t}, \gamma_{s,t}]=0.
\end{equation}
\end{proposition}

\begin{proof}
Note first of all that the right hand side of \eqref{eq:sum} defines a
Lie series, since only indices $k\le n$ contribute to the term of
degree $n$ in $\beta_{s,t}$. Also, $\beta_{0,t}=0$ since each term on 
the right hand side of \eqref{eq:sum} vanishes at $s=0$. 
The flatness condition \eqref{eq:flat} is satisfied as an equality of formal Lie series. To show that the series $\beta_{s,t}$
is strongly convergent, choose $D>0$ with $C_n( \gamma_{s,t}),\
C_n(\dot{\gamma}_{s,t})\le D^n$.  Using the estimate \eqref{eq:estimate}, we
obtain 
\[ \begin{split}
C_n(\beta_{s,t})&\le \sum_{k=0}^\infty 
\int_{0\leq s_0\leq s_1 \dots \leq s_k \leq s} ds_0 \dots ds_k\ 
C_n\big(\ad_{\gamma_{s_k,t}} \cdots \ad_{\gamma_{s_1,t}} \f{\p
  \gamma_{s_0,t} }{\p t}\big)\\
&\le \sum_{k=0}^\infty \f{s^k}{k!}  \f{2^k n^{2k}}{(2k-1)!!} \,D^n\le \sum_{k=0}^\infty \f{2^{2k} n^{2k}}{(2k)!} D^n 
\le e^{2n} D^n.\\
\end{split}
\]
The scaling property $\beta_{s,t}=t^{-1}r_t^*\beta_s$ holds since it 
is satisfied by each term in the sum \eqref{eq:sum}. 
\end{proof}

We now apply these results to the Torossian function, $\gamma_s=\ti{\beta}_s$.
The diagrammatic technique from Torossian's paper \cite{to:su} gives
$\ti{\beta}_s$ directly as a universal Lie series.  (The
family of maps $\ti{\beta}_s$ is however not canonical, since it depends
on the choice of a suitable path in $\ol{C}_{2,0}^+$, the compactified
configuration space of two points in the upper half plane.) The
estimates of \cite{an:ko} (see also Remark after Theorem 4.2 in
\cite{to:su}) show that this Lie series is strongly convergent (with
coefficients depending on the parameter $s$), and with $C(\ti{\beta}_s)$
uniformly bounded for $s\in [0,1]$. Define $\beta_{s,t}$ by the Lie
series \eqref{eq:sum}, and identify $\beta_{s,t}$ and $\ti{\beta}_{s,t}$
with the corresponding analytic functions on $\g\times\g$ (defined
near the origin). By the uniqueness property for solutions of the
zero curvature equation (cf. Lemma \ref{lem:zero}), the 
function $\beta_{1,t}=t^{-1}r_t^*\beta$ coincides with the KV solution
from Theorem \ref{th:suppl}. We have hence shown that our KV solution
is given by a universal, strongly convergent Lie series. 

\begin{remark}
One can replace any solution $\beta=(\beta^1,\beta^2)$ of the KV
conjecture with a {\em symmetric} solution
$\beta_{\on{sym}}=(\beta^1_{\on{sym}},\beta^2_{\on{sym}})$
where
\[
\begin{split}
\beta_{\on{sym}}^1(x,y)&=\hh \big(\beta^1(x,y)+\beta^2(-y,-x)\big),\\ 
\beta_{\on{sym}}^2(x,y)&=\hh \big(\beta^2(x,y)+\beta^1(-y,-x) \big) .
\end{split}
\]
This solution has the additional property $\beta^1_{\on{sym}}(x,y)=
\beta^2_{\on{sym}}(-y,-x)$.
\end{remark}

While Theorem \ref{th:suppl} settles the KV conjecture
in the form stated in \cite[page 250]{ka:ca}, there remain some
interesting open problems in this context. 
\begin{enumerate}
\item
The algebraic version III of the KV conjecture can also be considered
for Lie algebras over $\QQ$. We do not know, however, whether the
coefficients of our solution $\beta(x,y)$ are rational, for a suitable
choice of Torossian's function $\ti{\beta}_s(x,y)$.
%
%
\item
 In a recent paper, Alekseev and Petracci \cite{al:un} proved that for
universal symmetric solutions $\beta$ of the KV problem (for arbitrary
finite-dimensional Lie algebras over $\R$), the linear terms in $x$ or
in $y$ (i.e. the terms $(\partial_x\beta)(0,y)$ and
$(\partial_y\beta)(x,0)$) are uniquely determined. The linear terms in
Kashiwara and Vergne's solution for the case of solvable Lie algebras
\cite[Equation (5.2)]{ka:ca} are as prescribed by this result.  It is
unknown whether their formula might in fact give a solution for all
Lie algebras.
\end{enumerate}

\section{Lie algebra cohomology and the KV conjecture}
\label{sec:liecoh}
In this Section, we will give yet another version of the
Kashiwara-Vergne conjecture. This reformulation depends on a general
fact in Lie algebra cohomology, stated in Proposition
\ref{prop:trivial} below.
\subsection{Lie algebra cohomology}
For any $x\in\g$, 
we denote by $\iota^\wedge(x)$ the derivation of $\wedge\g^*$
given by contraction.  For $\mu\in\g^*$ we denote by $\eps(\mu)$ 
the operator of exterior multiplication on $\wedge\g^*$. 

For any $\g$-module $\M$, let $L^\M(x)\in \End(\M)$ denote the action
of $x\in \g$.  The Lie algebra cohomology $H^\bullet(\g,\M)$ of $\g$
with coefficients in $\M$ is the cohomology of the Chevalley-Eilenberg
complex
\begin{equation} C^k(\g,\M):=\M\otimes \wedge^k\g^*,\ \ 
\d=1\otimes\d^\wedge+\sum_a L^\M(e_a)\otimes \eps(e^a).\end{equation}
Here $e_a\in\g$ and $e^a\in\g^*$ are dual bases, and $\d^\wedge=\hh
\sum_a \eps(e^a) L^\wedge(e_a)$ is the Lie algebra differential on
$\wedge\g^*$. Suppose $\k$ is a another Lie algebra, with basis $f_i$, 
and that $\N$ is a $\k$-module.
Consider the linear map
\begin{equation} V\colon\End(\N)\otimes\k\to \End(\N),\ 
R=\sum_i R^i\otimes f_i \mapsto \sum_i R^i\circ
L^\N(f_i).\end{equation}
If $\psi\colon \g\to \k$ is a Lie algebra homomorphism and 
$R\in \on{End}(\N)\otimes \k$ is $\g$-invariant, then $V(R)\in
\End(\N)$ is $\g$-equivariant. Hence it defines a chain map, 
$V(R)\otimes\psi^*\colon C(\k,\N)\to C(\g,\N)$. 
We will need the following fact: 
\begin{proposition}\label{prop:trivial}
Let $\N$ be a $\k$-module, and 
$\psi\colon \g\to \k$ a Lie algebra homomorphism. Suppose
$R\in \on{End}(\N)\otimes \k$ is $\g$-invariant. Then the chain map
$V(R)\otimes \psi^*$ is homotopic to the trivial map. In fact,
\begin{equation}
V(R)\otimes \psi^*=(1\otimes\psi^*)\circ [\d,\iota(R)]
\end{equation}
with $\iota(R)=\sum_i R^i\otimes \iota(f_i)$. 
\end{proposition}

\begin{proof}
We have to show that $1\otimes \psi^*\colon \N\otimes\wedge\k^*\to
\N\otimes\wedge\g^*$ 
annihilates the expression, 
$[\d,\iota(R)]-V(R)\otimes 1$. Let $f^i\in\k^*$ be the dual basis to $f_i\in\k$.
We compute, 
\beq [\d,\iota(R)]-V(R)\otimes 1&=&[1\otimes \d^{\wedge}+\sum_j L^\N(f_j)\otimes \eps(f^j),
\sum_i R^i \otimes \iota(f_i)]-V(R)\otimes 1\\&=&\sum_i R^i \otimes L^\wedge(f_i)
+\sum_{ij}L^\N(f_j)(R^i)\otimes \big(\eps(f^j)\circ \iota(f_i)\big).
\eeq
Since $\sum_j (\psi^* f^j)\otimes f_j=\sum_a e^a\otimes \psi(e_a)$, 
we have the following equalities of maps $\wedge\k^*\to \wedge\g^*$,  
\beq 
\psi^*\circ L^\wedge(f_i)&=&
\sum_j \psi^* \circ \eps(f^j)\circ \iota([f_j,f_i]_\k)\\
&=&
\sum_j \eps(\psi^* f^j)\circ \psi^*\circ \iota([f_j,f_i]_\k)\\
&=& \sum_a \eps(e^a) \circ \psi^*\circ \iota([\psi(e_a),f_i]_\k).
\eeq
One the other hand, using the $\g$-equivariance of $R$,
\beq 
(1\otimes \psi^*)\circ 
\Big(\sum_{ij}L^\N(f_j)(R^i)\otimes \eps(f^j) \circ \iota(f_i)\Big)
&=& 
\sum_{ia}L^\N(\psi(e_a))(R^i)\otimes \eps(e^a)\circ  \psi^*\circ \iota(f_i)\\
&=&
-\sum_{ia}R^i\otimes  \big(\eps(e^a)\circ \psi^* \circ\iota([\psi(e_a),f_i]_\k)\big).
\eeq
\end{proof}

\subsection{Reformulation of the KV conjecture}\label{subsec:KV4}
We will apply this Proposition to the diagonal 
embedding $\psi\colon \g\to \g^2$. 
The map $\psi^*\colon \wedge\g^*\otimes\wedge\g^*=
\wedge(\g^*\oplus\g^*)\to \wedge\g^*$ is just the 
product map in the exterior algebra. View
$\N=\D'_{\on{comp}}(O^2)$ as a $\g^2$-module, and 
$\M=\D'_{\on{comp}}(\g)$ as a $\g$-module. We obtain a family of 
chain maps, 
\begin{equation}\label{eq:product0}
M_t=m_t\otimes \psi^*\colon C(\g^2,\D'_{\on{comp}}(O^2))\to 
C(\g,\D'_{\on{comp}}(\g)).
\end{equation}
Version I of the Kashiwara-Vergne conjecture says 
\begin{equation}\f{\p m_t}{\p t}=-m_t\circ V(\beta_t),\end{equation}
where $\beta_t\in C^\infty(M,\k)$. (The components $\beta_t^i$ are 
identified with the corresponding operators of multiplication on
$\D'_{\on{comp}}(O^2)$.) Tensoring with $\psi^*\colon \wedge^2\g^*\to
\wedge\g^*$, and using Proposition \ref{prop:trivial}, we obtain the
following equivalent reformulation of the KV conjecture: \\

\noindent {\bf Kashiwara-Vergne conjecture, version IV}: 
{\em Taking $O$ smaller if necessary, 
there exists a $\g$-equivariant 
function $\beta\colon O^2\to \g^2$, with $\beta(0,0)=0$, such that 
$\beta_t=t^{-1}r_t^*\beta$ satisfies 
\begin{equation}\label{eq:KV4}
\f{\p M_t}{\p t}=-M_t\circ [\d,\iota(\beta_t)].
\end{equation}
}\\
The algebra isomorphism $H(\g,S\g)\to H(\g,U\g)$
is an immediate consequence of this conjecture.
More generally, one obtains a similar statement for convolution of
germs of distributions. Let $A_1,A_2\subset \g$ be closed cones 
with $A_1\cap -A_2=\{0\}$, and consider the diagram 
\[ \begin{CD}
C(\g,\mf{D}'(\g,A_1))\times C(\g,\mf{D}'(\g,A_1)) @>>>
C(\g,\mf{D}'(\g,A_1+A_2))\\
@VVV @VVV\\
C(\g,\mf{D}'(G,A_1))\times C(\g,\mf{D}'(G,A_1)) @>>>
C(\g,\mf{D}'(G,A_1+A_2))
\end{CD}\]
where the horizontal maps are product maps, and the 
vertical maps are induced by the Duflo maps. Version IV of the Kashiwara-Vergne conjecture implies
that this diagram commutes up to a chain homotopy. In particular,
passing to cohomology one obtains a commutative diagram

\[ \begin{CD}  H(\g,\mf{D}'(\g,A_1))\otimes H(\g,\mf{D}'(\g,A_2))
@>>> H(\g,\mf{D}'(\g,A_1+A_2))\\ @VVV @VVV\\
H(\g,\mf{D}'(G,A_1))\otimes H(\g,\mf{D}'(G,A_2))@>>>
H(\g,\mf{D}'(G,A_1+A_2)) \\
\end{CD}
\]

\subsection{Higher homotopies}
Version IV of the KV conjecture shows that each map ${\rm Duf}_t
\otimes {\rm id}: S(\g) \otimes \wedge \g^* \to U(\g) \otimes \wedge
\g^*$ is an algebra isomorphism up to homotopy. Solutions of the KV
problem provide infinitesimal homotopies for family of
products $M_t$. One can ask whether in fact ${\rm Duf}_t \otimes
{\rm id}$ extend to an $A_\infty$-morphisms (for a definition see
\cite{ma:ho}) and to construct infinitesimal homotopies for its higher
components. Indeed, let 
\begin{equation} M_t^{(3)}=M_t\circ (M_t\otimes 1)\colon C(\g^3,\D'_{\on{comp}}(O^3))\to
C(\g,\D_{\on{comp}}(\g))\end{equation}
denote the triple product map (defined for $O$ sufficiently
small). Let $\beta_t$ denote a solution of the KV-problem. Taking the
$t$-derivative of $M_t^{(3)}=M_t\circ (M_t\otimes 1)$, we obtain
$\dot{M}_t^{(3)}=-[\d,h_t^{(3)}]$ where $h_t^{(3)}$ is the homotopy operator
\begin{equation}
h_t=M_t^{(3)}\circ (\iota(\beta_t)\otimes 1)+M_t\circ
\iota(\beta_t)\circ (M_t\otimes 1).
\end{equation}
On the other hand, starting with $M_t^{(3)}=M_t\circ (1\otimes M_t)$,
we get another homotopy operator 
\begin{equation}
\ti{h}_t=M_t^{(3)}\circ (1\otimes \iota(\beta_t))+M_t\circ
\iota(\beta_t)\circ (1\otimes M_t).
\end{equation}
The differences $\ti{h}_t-h_t$ are cochain maps of degree $-1$. In the
spirit of version $IV$ of the KV conjecture, one can formulate the
following new conjecture: {\em For $O \subset \g$ a sufficiently small neighborhood
of the origin, there exists a $\g$-equivariant map
$\beta^{(3)}\colon O^3 \to \wedge^2 \g^3$, with $\beta^{(3)}(0)=0$, such that
$\beta^{(3)}_t=t^{-1} r_t^* \beta^{(3)}$ satisfies
\begin{equation}
\ti{h}_t-h_t= - M_t^{(3)} \circ [d, \iota(\beta^{(3)}_t)] .
\end{equation}
}
In a similar fashion, one can conjecture the existence
of higher homotopies given by equivariant functions 
$\beta^{(n)}\colon O^n \to \wedge^{n-1} \g^n$.

\begin{appendix}
\section{The zero curvature equation}
In this appendix we prove Lemma \ref{lem:zero}, the zero curvature
equation for the Lie algebra $\k_M$. Let us first recall the more
general zero curvature equation for vector fields on a manifold $M$.
\begin{lemma}\label{lem:zero1}
Let $Y_{s,t}\in\mf{X}(M)$ a given 2-parameter family of vector fields,
depending smoothly on $s,t\in [0,1]$, and assume that $Y_{s,t}$
vanishes at some given point $x_0\in M$ for all $s,t$. Then there exists an
open neighborhood $U$ of $x_0$, and a unique 2-parameter family of
vector fields $X_{s,t}$, $s,t\in [0,1]$, satisfying the zero curvature
equation
\begin{equation} \label{eq:zerocurv}
\f{\p X_{s,t}}{\p s}-\f{\p Y_{s,t}}{\p t}+[X_{s,t},Y_{s,t}]=0
\end{equation}
with initial condition $X_{0,t}=0$. 
\end{lemma} 
\begin{proof}
Let $F_{s,t}$ denote the flow of
$Y_{s,t}$, viewed as an $s$-dependent vector field depending on $t$ as
a parameter. In terms of the action on smooth functions $f$,
\begin{equation}\label{eq:defy}
 (Y_{s,t}f)\circ F_{s,t}^{-1}=\f{\p}{\p s} (f\circ F_{s,t}^{-1});\ \ F_{0,t}=\on{id}.
\end{equation}
Since $Y_{s,t}(x_0)=0$, there exists an open neighborhood $U$ of $x_0$
such that $F_{s,t}(x)$ is defined for all $x\in U$ and $s,t\in [0,1]$.
Define a vector field $X_{s,t}$ by 
\begin{equation}
\label{eq:defx}
(X_{s,t}f)\circ F_{s,t}^{-1}=\f{\p}{\p t} (f\circ F_{s,t}^{-1}).
\end{equation}
Note that $X_{0,t}=0$ since $F_{0,t}=\on{id}$.
Subtracting the $s$-derivative of \eqref{eq:defx} from the
$t$-derivative of \eqref{eq:defy} one finds that $X_{s,t}$ solves 
\eqref{eq:zerocurv}. On the other hand, any two solutions of 
\eqref{eq:zerocurv} differ by a 2-parameter family of vector fields 
$Z_{s,t}$ with $\f{\p}{\p s} Z_{s,t}+[Z_{s,t},Y_{s,t}]=0$, or 
equivalently $\f{\p}{\p s}\big((F_{s,t}^{-1})^* Z_{s,t}\big)=0$. 
Integrating with initial condition $Z_{0,t}=0$ one obtains
$Z_{s,t}=0$, proving uniqueness. 
\end{proof}

Lemma \ref{lem:zero} (the zero curvature equation in $\k_M$)
may be viewed as a special case of Lemma \ref{lem:zero1}. Indeed, $\k_M$
may be realized as a Lie algebra of vector fields on $K\times M$,
where $K$ is a Lie group with Lie algebra $\k$. To see this, 
lift the Lie algebra $\k$-action on $M$ to an action on $K\times M$ by 
$\xi\mapsto \hat{X}^\xi=(\xi^L,X^\xi)$
where $\xi^L\in\mf{X}(K)$ is the left-invariant vector field generated by $\xi$. 
Accordingly, the homomorphism $\k_M\to \mf{X}(M)$ lifts to a Lie
algebra homomorphisms 
\begin{equation}\label{eq:inclusion}
\k_M\to \mf{X}(K\times M),\ \beta\mapsto
\hat{X}^\beta.
\end{equation} 
It is clear that the map \eqref{eq:inclusion} is injective. Its image
consists of vector fields on $K\times M$ that are invariant under 
the $K$-action by left-multiplication on the first factor, and are
tangent to the foliation defined by the $\k$-action $\xi\mapsto \hat{X}^\xi$
on $K\times M$. (Equivalently, the flow of such vector fields is 
$K$-equivariant, and preserves the leaves of the $\k$-action.) 
Using \eqref{eq:inclusion} to identify elements of $\k_M$ with vector
fields,  Lemma \ref{lem:zero} is now a direct consequence of 
\eqref{eq:zerocurv}. 
\end{appendix}

\bibliographystyle{amsplain}   %
\def\polhk#1{\setbox0=\hbox{#1}{\ooalign{\hidewidth
  \lowre1.5ex\hbox{`}\hidewidth\crcr\unhbox0}}} \def\cprime{$'$}
  \def\cprime{$'$} \def\polhk#1{\setbox0=\hbox{#1}{\ooalign{\hidewidth
  \lower1.5ex\hbox{`}\hidewidth\crcr\unhbox0}}} \def\cprime{$'$}
  \def\cprime{$'$}
\providecommand{\bysame}{\leavevmode\hbox to3em{\hrulefill}\thinspace}
\providecommand{\MR}{\relax\ifhmode\unskip\space\fi MR }
\providecommand{\MRhref}[2]{%
  \href{http://www.ams.org/mathscinet-getitem?mr=#1}{#2}
}
\providecommand{\href}[2]{#2}

\end{document}